\documentclass[10pt,a4paper]{article}

\usepackage[latin2]{inputenc}
\usepackage{amsmath}
\usepackage{amsfonts}
\usepackage{amssymb}
\usepackage{makeidx}

\begin{document}

$\,$

\bigskip

\begin{center}
\begin{Large}
\textbf{Reproductive and non-reproductive solutions \\[0.2 ex]
of the matrix equation {\boldmath $AXB=C$}}
\end{Large}
\end{center}

\smallskip
\begin{center}
{\sc Branko Male\v sevi\' c}${}^{\ast} $ and {\sc Biljana Radi\v ci\' c}
\end{center}

\smallskip
\begin{center}
{\small $^{\ast}$University of Belgrade, Faculty of Electrical Engineering, \\
Department of Applied Mathematics, Serbia}
\end{center}

\medskip
{\small \textbf{Abstract.}
In this article we consider a consistent matrix equation \mbox{$AXB\!=\!C$} when a particular
solution $X_{0}$ is given and we represent a new form of the general solution which contains
both reproductive and non-reproductive solutions (\mbox{it~depends} on the form of particular
solution $X_{0}$). We also analyse the solutions of some matrix systems using the concept of
reproductivity and we give a new form of the condition for the consistency of the matrix
equation $AXB=C$.}


\section{Reproductive equations}

The concept of the reproductive equations was introduced by {\sc S.$\,$B.$\;$Pre\v si\'c}~\cite{Presic68}.

\medskip
\textbf{Definition 1.2.} \textit{The reproductive equations} are the equations of the
following form:
\begin{equation}
x=f(x),
\end{equation}
where $x$ is a unknown, $S$ is a given set and $f:S \longrightarrow S $
is a given function which satisfies the following condition:
\begin{equation}
\label{UR}
f\circ f=f.
\end{equation}

\medskip
The condition (\ref{UR}) is called \textit{the condition of reproductivity}. The most
important statements in relation to the reproductive equations are given by the following
two theorems \cite{Presic68} (see also \cite{Presic72}, \cite{Bozic75} and \cite{Presic00}):

\medskip
\textbf{Theorem 1.1.} ({\sc S.$\,$B.$\;$Pre\v si\'c})
For any consistent equation $J(x)$ there is an equation of the form
$x=f(x)$,
which is equivalent to
$J(x)$
being in the same time reproductive as well. $\blacklozenge$

\smallskip
\textbf{Theorem 1.2.} ({\sc S.$\,$B.$\;$Pre\v si\'c})
If a certain equation $J(x)$ is equivalent to the reproductive one $x=f(x)$, the general solution
is given by the formula $x=f(y)$, for any value $y\in S$. $\blacklozenge$

\medskip
The concept of the reproductive equations allows us to analyse the solutions of some matrix
systems (see Application 2.1. and Application 2.2. in the following section of this paper). In
\cite{Keckic82}, \cite{Keckic83} and \cite{KeckicPresic97} authors considered the general
applications of the concept of reproductivity.

\section{The matrix equation {\boldmath $ AXB=C $}}

Let $m,n \in \mathbb{N}$ and $\mathbb{C}$ is the field of complex numbers. The set of all
matrices of order $m \times n$ over $\mathbb{C}$ is denoted by $\mathbb{C}^{m \times n}$.
For the set of all matrices from $\mathbb{C}^{m \times n}$ with a rank $a$ we use denotement
$\mathbb{C}_{a}^{m \times n}$ . Let $A=[a_{i,j}]\in \mathbb{C}^{m \times n}$.
By $A_{i\rightarrow}$ we denote the $i$-th row of $A$, $i=1,...,m$. For the $j$-th column of $A$,
$j=1,...,n$, we use denotement $A_{\downarrow j}$.

\medskip
A solution of the matrix equation
\begin{equation}
\label{AXAA}
AXA=A
\end{equation}
is called $\lbrace1\rbrace$-inverse of the matrix $A$ and it is denoted by $A^{(1)}$. The
set of all $\lbrace 1\rbrace$-inverses of the matrix $A$ is denoted by $ A\lbrace1\rbrace $.
For the matrix $A$, let regular matrices $Q\in \mathbb{C}^{m \times m}$ and
$P\in \mathbb{C}^{n \times n}$ be determined so that the following
equality is true.
\begin{equation}
\label{EA}
QAP=E_{a}=
\left[
\begin{array}{c|c}
I_{a} & 0\\\hline
0 & 0
\end {array}
\right],
\end{equation}
where $a=rank(A)$. In \cite{Rohde64} C. Rohde showed that the general $\lbrace1\rbrace$-inverse $A^{(1)}$ can be represented in the following form:
\begin{equation}
\label{AJEDAN}
A^{(1)} = P
\left[
\begin{array}{c|c}
I_{a} & X_{1}\\\hline
X_{2} & X_{3}
\end {array}
\right]Q,
\end{equation}
where $X_{1}$, $X_{2}$ and $X_{3}$ are arbitrary matrices of suitable sizes.

\bigskip
Let $A\in \mathbb{C}^{m \times n}$, $B\in \mathbb{C}^{p \times q}$ and
$C\in \mathbb{C}^{m \times q}$. In the paper \cite{Penrose55} {\sc R.$\;$Penrose} proved the following
theorem related to the matrix equation \begin{equation}
\label{AXBC}
AXB=C.
\end{equation}

\textbf{Theorem 2.1.} ({\sc R.$\;$Penrose}) The matrix equation (\ref{AXBC}) is consistent iff for
some choice of $\lbrace 1\rbrace$-inverses $A^{(1)}$ and $B^{(1)}$ of the matrices $A$ and $B$
the condition
\begin{equation}
\label{UKABC}
AA^{(1)}CB^{(1)}B=C
\end{equation}is true. The general solution of the matrix equation (\ref{AXBC})
is given by the formula
\begin{equation}
\label{ORR}
X=f(Y)=A^{(1)}CB^{(1)}+Y-A^{(1)}AYBB^{(1)},
\end{equation}
where  $ Y $ is an arbitrary matrix of suitable size. $\blacklozenge$

\medskip
\noindent
If a particular solution $X_{0}$ of the matrix equation (\ref{AXBC}) is given, the formula of
general solution is given in the following theorem.

\textbf{Theorem 2.2.} If $X_{0}$ is a particular solution of the matrix equation (\ref{AXBC}),
then the general solution of the matrix equation (\ref{AXBC}) is given by the formula
\begin{equation}
\label{OR}
X=g(Y)=X_{0}+Y-A^{(1)}AYBB^{(1)},
\end{equation}
where $Y$ is an arbitrary matrix of suitable size. The function $g$ satisfies the condition of
reproductivity (\ref{UR}) iff $X_{0}=A^{(1)}CB^{(1)}$.

\medskip
\textbf{Proof.} See \cite{MalesevicRadicic} and \cite{RadicicMalesevic} (where different proofs
are given).~$\blacklozenge$

\medskip
The formula (\ref{OR}) contains both reproductive and non-reproductive solutions. It depends on
the form of particular solution $X_{0}$.

\medskip
\textbf{Remark 2.1.}
In the paper \cite{MalesevicRadicic} authors proved that there is a matrix equation (\ref{AXBC})
and a particular solution $X_{0}$ so that:
\begin{equation}
\label{X0}
X_{0}\neq A^{(1)}CB^{(1)},
\end{equation}
for any choice of $\lbrace 1\rbrace$-inverses $A^{(1)}$ and $B^{(1)}$. In that case the formula
(\ref{OR}) gives the general non-reproductive solution. Otherwise, the formula (\ref{OR}) gives
the general reproductive solution. $\blacklozenge $

\medskip
\textbf{Example 2.1.} Compared to \cite{MalesevicRadicic}, we give a simpler example of
the matrix equation (\ref{AXBC}) and a particular solution $X_{0}$ so that (\ref{X0}) is valid.
Let's consider the matrix equation:
\begin{equation}
\label{AXBC_Ex_2_1}
\left[
\begin{array}{cc}
1 & 2
\end{array}
\right]
X
\left[
\begin{array}{c}
1 \\
3
\end{array}
\right]
=
[12],
\end{equation}
with
\mbox{\small $
A
\!=\!
\left[\!
\begin{array}{cc}
1 & 2
\end{array}
\!\right]
$},
\mbox{\small $
B
\!=\!
\left[\!
\begin{array}{cc}
1 & 3
\end{array}
\!\right]^{T}
\!$},
\mbox{\small $
C
\!=\!
\left[\!
\begin{array}{c}
12
\end{array}
\!\right]
$}.
Then, there is a particular solution:
\begin{equation}
\label{AXBC_Ex_2_1_Part_Sol}
X_{0}
=
\left[
\begin{array}{cc}
84 &-24 \\
-36 & 12
\end{array}
\right]
\end{equation}
of the previous matrix equation. It is easy to show that
\mbox{\small $
A^{(1)}
\!=\!
\left[
\begin{matrix}
1\!-\!2a & a
\end{matrix}
\right]^{T}\,(a \!\in\! \mathbb{C})$},
\mbox{\small $
B^{(1)}
\!=\!
\left[
\begin{matrix}
1\!-3b & b
\end{matrix}
\right]\,(b \!\in\! \mathbb{C}).
$}
So,

\vspace*{-2.0 mm}

\begin{equation}
X_{1}
=
A^{(1)}CB^{(1)}
=
\left[
\begin{matrix}
12\!-\!24a\!-\!36b\!+\!72ab & 12b\!-\!24ab \\
12a\!-\!36ab                & 12ab
\end{matrix}
\right].
\end{equation}
From
$
AXB-C
=0
\Longleftrightarrow
x_{1,1} +3 x_{1,2}+2x_{2,1}+6x_{2,2}-12 =0
$,
where
\mbox{\small $
X = \left[\!
\begin{array}{c}
\mbox{\normalsize $x_{i,j}$}
\end{array}
\!\right]
$},
we get that the matrix of general solution is given by the following form:
\begin{equation}
X =
\left[
\begin{matrix}
12-3p-2q-6r & p\\
q & r
\end{matrix}
\right]\!,\,
(p,q,r\in \mathbb{C}).
\end{equation}
For $p\!=\!-24, \, q\!=\!-36$ and $r\!=\!12$, we get the particular solution $X_{0}$ of the
matrix equation~(\ref{AXBC_Ex_2_1}), but $X_{0} \!\neq\! X_{1} \!=\! A^{(1)}CB^{(1)}$ for any
choice of $\lbrace 1\rbrace-$inverses $A^{(1)}$~and~$B^{(1)}$,
because from $X_{0} \!=\! X_{1}$ we obtain the contradiction ($ab \!=\! 1$~and $a \!=\! 0,\,b
\!=\! 0$).~$\blacklozenge$

\medskip
In \cite{Presic63} {\sc S.$\,$B.$\;$Pre\v si\' c} analysed the matrix equation (\ref{AXAA}) and he proved
the following theorem${}^{\natural)}$\footnote{${}^{\natural)}$ \scriptsize $\!\!$with the first
appearances of non-reproductive solutions {\big (}see $(E_3)$-$(E_5)${\big )}}$\!$.

\medskip
\textbf{Theorem 2.3.}
For any square matrix $A\in \mathbb{C}^{n \times n}$ and any
general~\mbox{$\lbrace 1\rbrace$-inverse} $A^{(1)}$ the following equivalences are true:

\vspace*{-5.0 mm}

$$
\begin{array}{l}
\mbox{\small $(E_1)\;\;$ $\;\; AX = 0 \; \Longleftrightarrow \;
(\exists Y \!\in \! \mathbb{C}^{n \times n}) \; X = Y -A^{(1)}AY$},                                      \\[0.2 ex]
\mbox{\small $(E_2)\;$  $\;\; XA = 0 \; \Longleftrightarrow \;
(\exists Y \!\in \!\mathbb{C}^{n \times n}) \; X = Y - YAA^{(1)}$},                                     \\[0.2 ex]
\mbox{\small $(E_3)$   $\;\; AXA = A \; \Longleftrightarrow \;
(\exists Y \!\in \!\mathbb{C}^{n \times n}) \; X = A^{(1)} + Y -A^{(1)}AYAA^{(1)}$},                    \\[0.2 ex]
\mbox{\small $(E_4)\;$  $\;\; AX = A \; \Longleftrightarrow \;
(\exists Y \!\in \! \mathbb{C}^{n \times n}) \; X = I + Y -A^{(1)}AY$},                                  \\[0.2 ex]
\mbox{\small $(E_5)\;\;$ $\;\; XA = A \; \Longleftrightarrow \;
(\exists Y \!\in \!\mathbb{C}^{n \times n}) \; X = I + Y - YAA^{(1)}$}. \; \blacklozenge
\end{array}
$$

In the general case the general solutions $(E_3)$-$(E_5)$ of Theorem 2.3. do not directly appear
according to {\sc Penrose}'s theorem. In \cite{Haveric83} {\sc M.$\;$Haveri\'c} showed that we can
get {\sc Penrose}'s solutions from {\sc Pre\v si\' c}'s solutions. She proved the following statement.

\break

\medskip
\textbf{Theorem 2.4.}
For any square matrix $A\in \mathbb{C}^{n \times n}$ and any general $\lbrace 1\rbrace-$inverse $A^{(1)}$ the following equivalences are true.

\vspace*{-5.0 mm}

$$
\begin{array}{l}
\mbox{\small $(E_1)\;\;$ $\;\; AX = 0 \; \Longleftrightarrow \;
(\exists Y \!\in \! \mathbb{C}^{n \times n}) \; X = Y - A^{(1)}AY$},                                     \\[0.2 ex]
\mbox{\small $(E_2)\;$  $\;\; XA = 0 \; \Longleftrightarrow \;
(\exists Y \!\in \!\mathbb{C}^{n \times n}) \; X = Y - YAA^{(1)}$},                                     \\[0.2 ex]
\mbox{\small $(E^{'}_3)$   $\;\; AXA = A \; \Longleftrightarrow \;
(\exists Y \!\in\!\mathbb{C}^{n \times n}) \; X = A^{(1)}AA^{(1)} + Y - A^{(1)}AYAA^{(1)}$},           \\[0.2 ex]
\mbox{\small $(E^{'}_4)\;$  $\;\; AX = A \; \Longleftrightarrow \;
(\exists Y \!\in \! \mathbb{C}^{n \times n}) \; X = A^{(1)}A + Y - A^{(1)}AY$},                          \\[0.2 ex]
\mbox{\small $(E^{'}_5)\;\;$ $\;\; XA = A \; \Longleftrightarrow \;
(\exists Y \!\in \! \mathbb{C}^{n \times n}) \; X = AA^{(1)} + Y - YAA^{(1)}$}. \; \blacklozenge
\end{array}
$$
\quad\enskip
Let's note that the previous two theorems are special case of Theorem~2.2.

\smallskip
For a consistent matrix equation (\ref{AXBC}) the following equivalence is true:
\begin{equation}
AXB=C \;\;\Longleftrightarrow\;\; X=f(X)=X-A^{(1)}(AXB-C)B^{(1)}.
\end{equation}
Therefore, based on Theorem 1.2., we have a short proof of the generality of formula (7) in
Theorem 2.1. (see \cite{RadicicMalesevic}).

\medskip
In the following applications we analysed the solutions of some matrix systems using the concept
of reproductivity.

\medskip
\textbf{Application 2.1.} Let $A,$ $B,$ $D$ and $E$ be given complex matrices of suitable sizes.
If the following matrix system is consistent:
\begin{equation}
\label{S1}
AX=B
\quad \wedge \quad
XD=E,
\end{equation}
then the general solution is given by the formula (\cite{Ben-IsraelGreville03}, {\sc A.$\;$Ben-Israel} and {\sc T.$\,$N.$\,$E.$\;$Greville})
\begin{equation}
X=g(Y)=X_{0}+(I-A^{(1)}A)Y(I-DD^{(1)}),
\end{equation}
where $Y$ is an arbitrary matrix of suitable size.

\medskip
In \cite{RadicicMalesevic} authors proved that if the matrix system (\ref{S1}) is consistent, the
general reproductive solution is given by the formula
\begin{equation}
\label{ORRS1}
X=f(Y)=A^{(1)}B+ED^{(1)}-A^{(1)}AED^{(1)}+(I-A^{(1)}A)Y(I-DD^{(1)}),
\end{equation}
where $Y$ is an arbitrary matrix of suitable size. The proof is based on the equivalence
\begin{equation}
\label{App_1}
(AX=B \;\wedge\; XD=E)
\;\;\Longleftrightarrow\;\;
X=f(X)
\end{equation} and Theorem 1.2. A more detailed proof can be found in \cite{RadicicMalesevic}.
$\blacklozenge$

\medskip
\textbf{Application 2.2.} Let $A \in \mathbb{C}^{n \times n}$ be a singular matrix.
If the following matrix system is consistent:
\begin{equation}
\label{Sys2}
AXA=A \quad \wedge \quad AX=XA,
\end{equation}
then the general solution is given by the formula (\cite{Keckic85}, {\sc J.$\,$D.$\;$Ke\v cki\' c})
\begin{equation}
\label{fSolA2}
X
=
f(Y)
=
Y + \bar{A} A \bar{A} - \bar{A} A Y - Y A \bar{A} + \bar{A} A Y A \bar{A},
\end{equation}
where $Y$ is an arbitrary matrix of suitable size and $\bar{A}$ is a commutative \mbox{$\{\!1\!\}$-inverse}.

\smallskip
In \cite{RadicicMalesevic} authors proved that (\ref{fSolA2}) represents the general solution of (\ref{Sys2}) using the concept of reproductivity. Further applications of the concept of reproductivity for some matrix equations and systems is considered in paper \cite{MalesevicRadicic12}.  $\blacklozenge$

\break

\medskip
The following theorem gives the new form of the condition for the consistency of the matrix equation (\ref{AXBC}). In the formulation of the theorem we use the following matrices
\begin{equation}
\label{ABCKAPA}
\widehat{A}=T_{A}A,\quad
\widehat{B}=BT_{B} \quad
\mbox{and} \quad
\widehat{C}=T_{A}CT_{B}
\end{equation}where $T_{A}$ is a permutation matrix which permutes linearly independent rows of
the matrix $A$ at the first $a$ positions and $T_{B}$ is a permutation matrix which
permutes linearly independent columns of the matrix $B$ at the first $b$ positions.

\medskip
Therefore, the matrix $\widehat{A}$ has linearly independent rows at the first $a$ positions and
the matrix $\widehat{B}$ has linearly independent columns at the first $b$ positions.

\smallskip
\noindent
Let
\begin{equation}
\label{VRSTEAKAPA}
\widehat{A}_{i \rightarrow}= \sum_{l=1}^{a}\alpha_{i,l} \widehat{A}_{l
\rightarrow},\quad i=a+1,...,m
\end{equation}
and
\begin{equation}
\label{KOLONEBKAPA}
\widehat{B}_{\downarrow j}=\sum_{k=1}^{b}\beta_{k,j}\widehat{B}_{\downarrow k},\qquad
j=b+1,...,q.
\end{equation}
for  some  scalars $\alpha_{i,l}$ and $\beta_{k,j}$. Then, the following theorem is true.

\medskip
\textbf{Theorem 2.5.} Let $A\!\in \!\mathbb{C}_{a}^{m \times n}$,
$B\!\in\!\mathbb{C}_{b}^{p \times q} $ and $C\!\in \!\mathbb{C}^{m \times q}$. Suppose that
$\widehat{A}$, $\widehat{B}$ and $\widehat{C}$ are determined by (\ref{ABCKAPA}) and that
(\ref{VRSTEAKAPA}) and (\ref{KOLONEBKAPA}) are satisfied. Then, the condition (\ref{UKABC})
is true for any choice of $\lbrace 1\rbrace$-inverses $A^{(1)}$ and $B^{(1)}$ iff

\vspace*{3 mm}

\noindent
$\!\!\!\!\widehat{C}$=$\left[
\begin{array}{cccccc}
c_{1,1}\!&\!...\!&\!c_{1,b}\!&\!\mbox{\tiny
$\displaystyle\sum_{k=1}^{b}$}\beta_{k,b+1}c_{1,k}\!&\!...\!&\mbox{\tiny
$\displaystyle\sum_{k=1}^{b}$}\beta_{k,q}c_{1,k}\\[-0.2 ex]
 .\!&\!...\!&.&.&\!...\!&.\\[-0.2 ex]
 .\!&\!...\!&.&.&\!...\!&.\\[-0.2 ex]
 .\!&\!...\!&.&.&\!...\!&.\\[-0.2 ex]
 c_{a,1}\!&\!...\!&c_{a,b}&\mbox{\tiny $\displaystyle\sum_{k=1}^{b}$}\beta_{k,b+1}c_{a,k}&\!...
\!&\mbox{\tiny $\displaystyle\sum_{k=1}^{b}$}\beta_{k,q}c_{a,k} \\[-0.2 ex]
\mbox{\tiny $\displaystyle\sum_{l=1}^{a}$}\alpha_{a+1,l}c_{l,1}\!&\!...\!& \mbox{\tiny
$\displaystyle\sum_{l=1}^{a}$}\alpha_{a+1,l}c_{l,b}\!&
\begin{small}\mbox{\tiny $\displaystyle\sum_{l=1}^{a}$}\mbox{\tiny
$\displaystyle\sum_{k=1}^{b}$}\alpha_{a+1,l}\beta_{k,b+1}c_{l,k}\end{small} &\!...\!&
\begin{small}\mbox{\tiny $\displaystyle\sum_{l=1}^{a}$}\mbox{\tiny $\displaystyle\sum_{k=1}^{b}$}
\alpha_{a+1,l}\beta_{k,q}c_{l,k}\end{small}\\[-0.2 ex]
 .&\!...\!& .& .&\!...\!& .\\[-0.2 ex]
 .&\!...\!& .& .&\!...\!& .\\[-0.2 ex]
 .&\!...\!& .& .&\!...\!& .\\[-0.2 ex]
\mbox{\tiny $\displaystyle\sum_{l=1}^{a}$}\alpha_{m,l}c_{l,1}&\!...\!& \mbox{\tiny
$\displaystyle\sum_{k=1}^{b}$}\alpha_{m,l}c_{l,b}&
\begin{small}\mbox{\tiny $\displaystyle\sum_{l=1}^{a}$}\mbox{\tiny
$\displaystyle\sum_{k=1}^{b}$}\alpha_{m,l}\beta_{k,b+1}c_{l,k} \end{small}&\!...\!&
\begin{small} \mbox{\tiny $\displaystyle\sum_{l=1}^{a}$}\mbox{\tiny
$\displaystyle\sum_{k=1}^{b}$}\alpha_{m,l}\beta_{k,q}c_{l,k} \end{small}
\end {array}
\right]$,

\vspace*{3 mm}

\noindent
where $c_{i,j}$ are arbitrary elements of $\mathbb{C}$.

\vspace*{3 mm}

\textbf{Proof.} The proof can be found in \cite{RadicicMalesevic}. $\blacklozenge$

\vspace*{3 mm}

\noindent
\textbf{Acknowledgements.}
The research is partially supported by the Ministry of Education and Science, Serbia, Grant No.
174032.


\begin{thebibliography}{00}

\bibitem{Penrose55} {\sc R.$\;$Penrose}: {\em A generalized inverses for matrices}, Math. Proc. Cambridge
Philos. Soc. 51 (1955), 406-413.

\vspace*{0.5 mm}

\bibitem{Presic63} {\sc S.$\,$B.$\;$Pre\v si\' c}: {\em Certaines  \' equations matricielles}, Publ. Elektrotehn.
Fak. Ser. Mat.-Fiz. $\mathit{N}^{\underline{o}}$ 121, Beograd, 1963. ({\tt http:/$\!$/pefmath.etf.rs/})

\vspace*{0.5 mm}

\bibitem{Rohde64} {\sc C.$\,$A.$\;$Rohde}: {\em Contribution to the theory, computation and application of
generalized inverses}, Doctoral  dissertation, University of North Carolina at Releigh, May 1964.

\vspace*{0.5 mm}

\bibitem{Presic68} {\sc S.$\,$B.$\;$Pre\v si\' c}: {\em Une classe d'\' equations matricielles et l'\' equation
fonctionnelle $f^{2}\!=\!f$}, Publications de l'institut mathematique, Nouvelle serie, T. 8
(22), Beograd 1968, 143-148. ({\tt http:/$\!$/publications.mi.sanu.ac.rs/})

\vspace*{0.5 mm}

\bibitem{Presic72} {\sc S.$\,$B.$\;$Pre\v si\' c}: {\em Ein Satz \" Uber Reproduktive L\" osungen},
Publications de l'in\-sti\-tut mathematique, Nouvelle serie, T. 14 (28), Beograd 1972, 133-136.

\vspace*{0.5 mm}

\bibitem{Bozic75} {\sc M.$\;$Bo\v zi\' c}: {\em A Note On Reproductive Solutions}, Publications  de l'institut
mathematique, Nouvelle serie, T. 19 (33), Beograd 1975,~33-35.

\vspace*{0.5 mm}

\bibitem{Keckic82} {\sc J.$\,$D.$\;$Ke\v cki\' c}: {\em Reproductivity of some equations of analysis {\rm I}},
Publications  de l'institut  mathematique, Nouvelle serie, T. 31 (45), Beograd 1982, 73-81.

\vspace*{0.5 mm}

\bibitem{Keckic83} {\sc J.$\,$D.$\;$Ke\v cki\' c}: {\em Reproductivity of some equations of analysis {\rm II}},
Publications  de l'institut  mathematique, Nouvelle serie, T. 33 (47), Beograd 1983, 109-118.

\vspace*{0.5 mm}

\bibitem{Haveric83} {\sc M.$\;$Haveri\' c}: {\em Formulae for general reproductive solutions of certain matrix
equations}, Publications de l'institut mathematique, Nouvelle serie, T. 34 (48), Beograd 1983,
81-84.

\vspace*{0.5 mm}

\bibitem{Keckic85} {\sc J.$\,$D.$\;$Ke\v cki\' c}: {\em Commutative weak generalized inverses of a square matrix
and  some  related  matrix   equations}, Publications  de l'institut  mathematique, Nouvelle
serie, T. 38 (52), Beograd 1985, 39-44.

\vspace*{0.5 mm}

\bibitem{KeckicPresic97} {\sc J.$\,$D.$\;$Ke\v cki\' c} and {\sc S.$\,$B.$\;$Pre\v si\' c}: {\em Reproductivity$\;-\,$A
general approach to equations}, Facta Universitatis (Ni\v s), Ser. Math. Inform. 12 (1997),
\mbox{157-184}.

\vspace*{0.5 mm}

\bibitem{Presic00} {\sc S.$\,$B.$\;$Pre\v si\' c}: {\em A generalization of the notion of reproductivity},
Publications de l'institut mathematique, Nouvelle serie, T. 67
(81), Beograd 2000, 76-84

\vspace*{0.5 mm}

\bibitem{Ben-IsraelGreville03} {\sc A.$\;$Ben-Israel} and {\sc T.$\,$N.$\,$E.$\;$Greville}: {\em Generalized Inverses$:$ Theory
and Applications}, Springer, 2003.

\vspace*{0.5 mm}

\bibitem{MalesevicRadicic} {\sc B.$\;$Male\v sevi\' c} and {\sc B.$\;$Radi\v ci\' c}: {\em Non-reproductive~and~\mbox{reproductive}~solu\-tions
of some matrix equations}, Proc. of Inter. Conf. "{\sl Mathematical and Infor\-mational Technologies}$\,$"$\!$, \mbox{\sl MIT-2011}$\,$,
{\sl http/$\!$/www.mit.rs/}, Vrnja\v cka Banja, Serbia, August 2011. ({\tt http:/$\!$/arxiv.org/abs/1106.3818})

\vspace*{0.5 mm}

\bibitem{MalesevicRadicic12} {\sc B. Male\v sevi\' c} and {\sc B. Radi\v ci\' c}, {\em Some considerations of matrix equations
using the concept of reproductivity}, appear in Kragujevac Journal of Mathematics.
({\tt http:/$\!$/arxiv.org/abs/1110.3519})

\vspace*{0.5 mm}

\bibitem{RadicicMalesevic} {\sc B.$\;$Radi\v ci\' c} and {\sc B.$\;$Male\v sevi\' c}: {\em Some considerations
in relation to the matrix equation $AXB\!=\!C$}, preprint. ({\tt http:/$\!$/arxiv.org/abs/1108.2485})

\end{thebibliography}
\end{document}